\documentclass{amsproc}%
\usepackage{graphicx}
\usepackage{amscd}
\usepackage{amsmath}
\usepackage{amsfonts}
\usepackage{amssymb}%
\theoremstyle{plain}

\numberwithin{equation}{section}

\begin{document}
\title[Matrix representations of finitely generated Grassmann algebras]{Matrix representations of finitely generated Grassmann algebras and some consequences}
\author{L\'{a}szl\'{o} M\'{a}rki}
\address{Alfr\'{e}d R\'{e}nyi Institute of Mathematics, Hungarian Academy of Sciences,
1364 Budapest, Pf. 127, Hungary}
\email{marki.laszlo@renyi.mta.hu}
\author{Johan Meyer}
\address{Department of Mathematics and Applied Mathematics, University of the Free
State, PO Box 339, Bloemfontein 9300, South Africa}
\email{MeyerJH@ufs.ac.za}
\author{Jen\H{o} Szigeti}
\address{Institute of Mathematics, University of Miskolc, 3515
Miskolc-Egyetemv\'{a}ros, Hungary}
\email{matjeno@uni-miskolc.hu}
\author{Leon van Wyk}
\address{Department of Mathematical Sciences, Stellenbosch University P/Bag
X1\noindent\noindent\noindent\noindent\ Matieland 7602, Stellenbosch, South
Africa }
\email{LvW@sun.ac.za}
\thanks{The first and third named authors were supported by OTKA K-101515 of Hungary}
\thanks{The second and fourth named authors were supported by the National Research
Foundation of South Africa under Grant No.~UID 72375. Any opinion, findings
and conclusions or recommendations expressed in this material are those of the
authors and therefore the National Research Foundation does not accept any
liability in regard thereto.}
\thanks{As far as the third author is concerned, this research was partly carried out
as part of the TAMOP-4.2.1.B-10/2/KONV-2010-0001 project with support from the
European Union, co-financed by the European Social Fund.}
\thanks{The authors thank P. N. \'Anh for fruitful discussions.}

\begin{abstract}
We prove that the $m$-generated Grassmann algebra can be embedded into a
$2^{m-1}\times2^{m-1}$ matrix algebra over a factor of a commutative
polynomial algebra in $m$ indeterminates. Cayley--Hamilton and standard
identities for $n\times n$ matrices over the $m$-generated Grassmann algebra
are derived from this embedding. Other related embedding results are also presented.

\end{abstract}
\subjclass{15A75, 16G30, 16R10}
\keywords{Grassmann algebra, matrix algebra, $K$-embedding, Cayley--Hamilton identity,
skew polynomial algebra}
\maketitle

\noindent1. INTRODUCTION

\bigskip

Let $K$ be a field (of characteristic zero) and consider the free associative
$K$-algebra $K\left\langle x_{1},\ldots,x_{m}\right\rangle $ generated by the
(non-commuting) indeterminates $x_{1},\ldots,x_{m}$. The elements of
$K\left\langle x_{1},\ldots,x_{m}\right\rangle $ are $K$-linear combinations
of monomials of the form $x_{i_{1}}\cdots x_{i_{k}}$ with $i_{1},\ldots
,i_{k}\in\{1,\ldots,m\}$ (not necessarily different). The $m$-generated
Grassmann (or exterior) algebra is defined as the factor $K$-algebra%
\[
E^{(m)}=K\left\langle x_{1},\ldots,x_{m}\right\rangle \diagup I(x_{1}%
,\ldots,x_{m}),
\]
where%
\[
I(x_{1},\ldots,x_{m})=(x_{i}x_{j}+x_{j}x_{i}\mid1\leq i\leq j\leq m)_{K}%
\]
is the (two-sided) $K$-ideal of $K\left\langle x_{1},\ldots,x_{m}\right\rangle
$\ generated by the polynomials $x_{i}x_{j}+x_{j}x_{i}$ for $1\leq i\leq j\leq
m$. The usual notation for this Grassmann algebra is%
\[
E^{(m)}=K\left\langle v_{1},\ldots,v_{m}\mid v_{i}v_{j}+v_{j}v_{i}=0\text{ for
all }1\leq i\leq j\leq m\right\rangle ,
\]
where the cosets $v_{i}=x_{i}+I(x_{1},\ldots,x_{m})$ are the anticommuting
generators of $E^{(m)}$. Clearly, any element of $E^{(m)}$ is a unique
$K$-linear combination of monomials of the form $v_{i_{1}}\cdots v_{i_{k}}$
with $1\leq i_{1}<\cdots<i_{k}\leq m$. It follows that $\dim_{K}E^{(m)}=2^{m}$.

The definition of the countably infinitely generated Grassmann algebra%
\[
E=K\left\langle v_{1},\ldots,v_{m},\ldots\mid v_{i}v_{j}+v_{j}v_{i}=0\text{
for all }1\leq i\leq j\right\rangle
\]
is similar.

For a ring (or $K$-algebra) $R$ let $\mathrm{M}_{n}(R)$ denote the full
$n\times n$ matrix ring ($K$-algebra) over $R$ with identity $I_{n}%
\in\mathrm{M}_{n}(R)$. Since any matrix $A\in\mathrm{M}_{n}(E)$ has a finite
number of entries and each entry contains a finite number of generators
$v_{i}$ from $E$, there exists an integer $m\geq1$ such that $A\in
\mathrm{M}_{n}(E^{(m)})$. It follows that%
\[
\mathrm{M}_{n}(E)=\overset{\infty}{%
{\textstyle\bigcup\limits_{m=1}}
}\mathrm{M}_{n}(E^{(m)}).
\]
The algebras $E^{(m)}$ and $E$ play a fundamental role in many areas of mathematics.

Our main inspiration is Kemer's pioneering work [K] on the $T$-ideals of
associative algebras, which revealed the importance of the identities
satisfied by the full $n\times n$ matrix algebra $\mathrm{M}_{n}(E)$ and by
the algebra $\mathrm{M}_{n,t}(E)$ of $(n,t)$-supermatrices (this is a certain
$K$-subalgebra of $\mathrm{M}_{n}(E)$). The prime $T$-ideals of $K\left\langle
x_{1},\ldots,x_{m},\ldots\right\rangle $\ are exactly the $T$-ideals of the
identities satisfied by $\mathrm{M}_{n}(K)$ for $n\geq1$. The $T$-prime
$T$-ideals are the prime $T$-ideals plus the $T$-ideals of the identities of
$\mathrm{M}_{n}(E)$ for $n\geq1$ and of $\mathrm{M}_{n,t}(E)$ for $n-1\geq
t\geq1$. Another remarkable result is that, for $n$ sufficiently large, any
$T$-ideal contains the $T$-ideal of the identities satisfied by $\mathrm{M}%
_{n}(E)$. Thus the algebras $\mathrm{M}_{n}(E)$ and $\mathrm{M}_{n}(E^{(m)})$
served as the main motivation for the present work.

An additional motivation is the well-known embedding of the skew field
$\mathbb{H}$\ of the real quaternions into $4\times4$ real matrices:%
\[
a+bi+cj+dk\,\,\longmapsto\left[
\begin{array}
[c]{cccc}%
\;\;a & \;\;b & \;\;c & \;\;d\\
-b & \;\;a & -d & \;\;c\\
-c & \;\;d & \;\;a & -b\\
-d & -c & \;\;b & \;\;a
\end{array}
\right]  ,
\]
where $a,b,c,d\in\mathbb{R}$. The above definition provides an injective
$\mathbb{R}$-algebra homomorphism $\upsilon:\mathbb{H}\rightarrow
\mathrm{M}_{4}(\mathbb{R})$. Using the natural extension%
\[
\upsilon_{n}:\mathrm{M}_{n}(\mathbb{H)}\longrightarrow\mathrm{M}%
_{n}(\mathrm{M}_{4}(\mathbb{R}))\cong\mathrm{M}_{4n}(\mathbb{R}),
\]
an $n\times n$ matrix over $\mathbb{H}$ can be viewed as a $4n\times4n$ matrix
over $\mathbb{R}$. For a quaternionic matrix $A\in\mathrm{M}_{n}(\mathbb{H})$,
the Cayley--Hamilton identity for $\upsilon_{n}(A)$ yields the same identity
(with real coefficients) of degree $4n$ for $A$ itself.

A similar approach to get a Cayley--Hamilton identity of degree $2n$\ for a
matrix $A\in\mathrm{M}_{n}(E^{(2)})$ is based on embedding the two-generated
exterior algebra $E^{(2)}$\ into a $2\times2$ matrix algebra over a certain
commutative ring (see [SzvW]).

In order to present a Cayley--Hamilton identity of degree $2^{m-1}n$ for a
matrix in $\mathrm{M}_{n}(E^{(m)})$, first we consider the so called constant
trace (CT-)representations of an arbitrary $K$-algebra. In Section 2 we derive
a Cayley--Hamilton identity for $n\times n$ matrices over any algebra having a
CT-representation. Then in Section 3, using induction, we present a
CT-representation ($K$-embedding)%
\[
\varepsilon^{(m)}:E^{(m)}\longrightarrow\mathrm{M}_{2^{m-1}}(K[z_{1}%
,\ldots,z_{m}]/(z_{1}^{2},\ldots,z_{m}^{2})).
\]
Notice that we cannot expect similar results for the infinitely generated
Grassmann algebra $E$. Since $E$ does not satisfy any of the standard
identities, it follows that $E$\ does not embed into any full matrix algebra
over a commutative ring.

Finally in Section 4 we use a certain factor of a skew polynomial algebra to
give a broad generalization of the embedding process for $E^{(m)}$.

\bigskip

\noindent2. CT-REPRESENTATIONS AND CAYLEY--HAMILTON IDENTITIES

\bigskip

Let $_{K}R$ be an arbitrary and $_{K}\Omega$ be a commutative (associative)
algebra over $K$ (notice that $K\subseteq\mathrm{Z}(R)$ and $K\subseteq\Omega
$). For an integer $t\geq1$, we consider representations of $R$ over $\Omega$
which are injective $K$-algebra homomorphisms ($K$-embeddings) $\varepsilon
:R\longrightarrow\mathrm{M}_{t}(\Omega)$. We call $\varepsilon$ a constant
trace (CT-) representation if $\mathrm{tr}(\varepsilon(r))\in K$ for all $r\in
R$ (here $\mathrm{tr}(\varepsilon(r))$ is the sum of the diagonal entries of
the $t\times t$ matrix $\varepsilon(r)\in\mathrm{M}_{t}(\Omega)$).

\bigskip

\noindent\textbf{2.1.Theorem.}\textit{ Let }$\varepsilon:R\longrightarrow
\mathrm{M}_{t}(\Omega)$\textit{ be a CT-representation of }$R$\textit{ over
}$\Omega$\textit{. If }$A\in\mathrm{M}_{n}(R)$\textit{ is an }$n\times
n$\textit{ matrix, then }$A$\textit{ satisfies a Cayley--Hamilton identity of
the form}%
\[
A^{tn}+c_{1}A^{tn-1}+\cdots+c_{tn-1}A+c_{tn}I_{n}=0,
\]
\textit{where }$c_{i}\in K$\textit{, }$1\leq i\leq tn$\textit{.}

\bigskip

\noindent\textbf{Proof.} Let%
\[
\varepsilon_{n}:\mathrm{M}_{n}(R)\longrightarrow\mathrm{M}_{n}(\mathrm{M}%
_{t}(\Omega))\cong\mathrm{M}_{tn}(\Omega)
\]
be the natural extension of $\varepsilon$. For any matrix $A=[a_{i,j}]$ in
$\mathrm{M}_{n}(R)$, the trace of the $tn\times tn$ matrix $B=\varepsilon
_{n}(A)$ is the sum of the traces of the diagonal $t\times t$ blocks:%
\[
\mathrm{tr}(B)=\overset{n}{\underset{i=1}{\sum}}\mathrm{tr}(\varepsilon
(a_{i,i})).
\]
Since $\varepsilon$ is a CT-representation, we have $\mathrm{tr}%
(\varepsilon(a_{i,i}))\in K$ for each $1\leq i\leq n$. It follows that
$\mathrm{tr}(B)\in K$. For the coefficients of the characteristic polynomial%
\[
\det(zI-B)=c_{0}z^{tn}+c_{1}z^{tn-1}+\cdots+c_{tn-1}z+c_{tn}\in\Omega\lbrack
z]
\]
of $B$, the following recursion holds: $c_{0}=1$ and%
\[
c_{k}=-\frac{1}{k}\left(  c_{k-1}\mathrm{tr}(B)+c_{k-2}\mathrm{tr}%
(B^{2})+\cdots+c_{1}\mathrm{tr}(B^{k-1})+c_{0}\mathrm{tr}(B^{k})\right)
\]
for $1\leq k\leq tn$ (Newton formulae, see [R]). In view of%
\[
\mathrm{tr}(B^{k})=\mathrm{tr}((\varepsilon_{n}(A))^{k})=\mathrm{tr}%
(\varepsilon_{n}(A^{k}))\in K,
\]
we deduce that $c_{i}\in K$ for each $0\leq i\leq tn$. Thus $\det(zI-B)\in
K[z]$ and the Cayley--Hamilton identity for $B\in\mathrm{M}_{tn}(\Omega)$ is
of the form%
\[
B^{tn}+c_{1}B^{tn-1}+\cdots+c_{tn-1}B+c_{tn}I_{n}=0.
\]
It follows that%
\[
\left(  \varepsilon_{n}(A)\right)  ^{tn}+c_{1}\left(  \varepsilon
_{n}(A)\right)  ^{tn-1}+\cdots+c_{tn-1}\varepsilon_{n}(A)+c_{tn}I_{n}=
\]%
\[
\varepsilon_{n}(A^{tn}+c_{1}A^{tn-1}+\cdots+c_{tn-1}A+c_{tn}I_{n})=0
\]
holds in $\mathrm{M}_{tn}(\Omega)$, thus the injectivity of $\varepsilon_{n}$
gives the desired identity. $\square$

\newpage

\noindent3. CT-REPRESENTATION\ OF\ $E^{(m)}$

\bigskip

For $m=1$ we have a natural isomorphism $E^{(1)}=K\left\langle v_{1}\mid
v_{1}^{2}=0\right\rangle \cong K[z]/(z^{2})$ of $K$-algebras. If $m=2$, then
the assignments%
\[
1\longmapsto\left[
\begin{array}
[c]{cc}%
1 & 0\\
0 & 1
\end{array}
\right]  ,\,\,v_{1}\longmapsto\left[
\begin{array}
[c]{cc}%
z_{1} & 0\\
0 & -z_{1}%
\end{array}
\right]  ,\,\,v_{2}\longmapsto\left[
\begin{array}
[c]{cc}%
0 & z_{2}\\
z_{2} & 0
\end{array}
\right]
\]
define the following CT-representation $\varepsilon^{(2)}:E^{(2)}%
\longrightarrow\mathrm{M}_{2}(K[z_{1},z_{2}]/(z_{1}^{2},z_{2}^{2}))$ :%
\[
\varepsilon^{(2)}(c_{0}\!+\!c_{1}v_{1}\!+\!c_{2}v_{2}\!+\!c_{3}v_{1}%
v_{2})\!=\!\left[
\begin{array}
[c]{cc}%
c_{0}+c_{1}z_{1}+(z_{1}^{2},z_{2}^{2}) & c_{2}z_{2}+c_{3}z_{1}z_{2}+(z_{1}%
^{2},z_{2}^{2})\\
c_{2}z_{2}-c_{3}z_{1}z_{2}+(z_{1}^{2},z_{2}^{2}) & c_{0}-c_{1}z_{1}+(z_{1}%
^{2},z_{2}^{2})
\end{array}
\right]  ,
\]
where $c_{0},c_{1},c_{2},c_{3}\in K$ and $(z_{1}^{2},z_{2}^{2})$ is the ideal
of the commutative polynomial ring $K[z_{1},z_{2}]$ generated by the monomials
$z_{1}^{2},z_{2}^{2}$.

\bigskip

\noindent\textbf{3.1.Theorem.}\textit{ For some integers }$m,t\geq2$, \textit{
let }$\varepsilon^{(m)}:E^{(m)}\longrightarrow\mathrm{M}_{t}(\Omega)$\textit{
be a CT-representation of }$E^{(m)}$\textit{ over a commutative }%
$K$\textit{-algebra }$\Omega$\textit{. Then the assignments}%
\[
1\longmapsto\left[
\begin{array}
[c]{cc}%
I_{t} & 0\\
0 & I_{t}%
\end{array}
\right]  ,\,\,v_{i}\longmapsto\left[
\begin{array}
[c]{cc}%
\varepsilon^{(m)}(v_{i}) & 0\\
0 & -\varepsilon^{(m)}(v_{i})
\end{array}
\right]  \text{\textit{ for }}1\leq i\leq m,
\]%
\[
\text{\textit{and }}v_{m+1}\longmapsto\left[
\begin{array}
[c]{cc}%
0 & \widehat{z}I_{t}\\
\widehat{z}I_{t} & 0
\end{array}
\right]  \text{\textit{ (with }}\widehat{z}=z+(z^{2})\text{\textit{ in }%
}\Omega\lbrack z]/(z^{2})\text{\textit{)}}%
\]
\textit{define a CT-representation }$\varepsilon^{(m+1)}:E^{(m+1)}%
\longrightarrow\mathrm{M}_{2t}(\Omega\lbrack z]/(z^{2}))$\textit{.}

\bigskip

\noindent\textbf{Proof.} Any element of $E^{(m+1)}$ can be uniquely written as
$g+hv_{m+1}$, where
\[
g=\underset{1\leq i_{1}<\ldots<i_{k}\leq m}{%
{\textstyle\sum}
}c_{i_{1},\ldots,i_{k}}v_{i_{1}}\cdots v_{i_{k}}\text{ and }h=\underset{1\leq
j_{1}<\ldots<j_{l}\leq m}{%
{\textstyle\sum}
}d_{j_{1},\ldots,j_{l}}v_{j_{1}}\cdots v_{j_{l}}%
\]
are in $E^{(m)}$\ with $c_{i_{1},\ldots,i_{k}},d_{j_{1},\ldots,j_{l}}\in K$.
For $1\leq i_{1}<\cdots<i_{k}\leq m$\ and $1\leq j_{1}<\cdots<j_{l}\leq m$,
our assignment gives%
\[
v_{i_{1}}\cdots v_{i_{k}}\longmapsto\left[
\begin{array}
[c]{cc}%
\varepsilon^{(m)}(v_{i_{1}}\cdots v_{i_{k}}) & 0\\
0 & (-1)^{k}\varepsilon^{(m)}(v_{i_{1}}\cdots v_{i_{k}})
\end{array}
\right]
\]
and%
\[
v_{j_{1}}\cdots v_{j_{l}}v_{m+1}\longmapsto\left[
\begin{array}
[c]{cc}%
0 & \varepsilon^{(m)}(v_{j_{1}}\cdots v_{j_{l}})\widehat{z}\\
(-1)^{l}\varepsilon^{(m)}(v_{j_{1}}\cdots v_{j_{l}})\widehat{z} & 0
\end{array}
\right]  .
\]
Thus%
\[
\varepsilon^{(m+1)}(g+hv_{m+1})=\left[
\begin{array}
[c]{cc}%
\varepsilon^{(m)}(g_{0}+g_{1}) & \varepsilon^{(m)}(h_{0}+h_{1})\widehat{z}\\
\varepsilon^{(m)}(h_{0}-h_{1})\widehat{z} & \varepsilon^{(m)}(g_{0}-g_{1})
\end{array}
\right]  ,
\]
where $g=g_{0}+g_{1}$ and $h=h_{0}+h_{1}$ are the unique presentations as sums
of an even and an odd element (with respect to the natural $\mathbb{Z}_{2}%
$-grading $E^{(m)}=E_{0}^{(m)}\oplus E_{1}^{(m)}$). Straightforward
verification shows that $\varepsilon^{(m+1)}:E^{(m+1)}\longrightarrow
\mathrm{M}_{2t}(\Omega\lbrack z]/(z^{2}))$ is an injective homomorphism of
$K$-algebras. In view of%
\[
\mathrm{tr}(\varepsilon^{(m+1)}(g+hv_{m+1}))=\mathrm{tr}(\varepsilon
^{(m)}(g_{0}+g_{1}))+\mathrm{tr}(\varepsilon^{(m)}(g_{0}-g_{1}))\in K,
\]
we deduce that $\varepsilon^{(m+1)}$ is a CT-representation of $E^{(m+1)}$.
$\square$

\bigskip

\noindent\textbf{3.2.Corollary.} \textit{For any integer }$m\geq2$\textit{
there exists a CT-representation }%
\[
\varepsilon^{(m)}:E^{(m)}\longrightarrow\mathrm{M}_{2^{m-1}}(K[z_{1}%
,\ldots,z_{m}]/(z_{1}^{2},\ldots,z_{m}^{2}))
\]
\textit{of }$E^{(m)}$\textit{, where }$(z_{1}^{2},\ldots,z_{m}^{2})$\textit{
is the ideal of }$K[z_{1},\ldots,z_{m}]$\textit{ generated by the monomials
}$z_{1}^{2},\ldots,z_{m}^{2}$\textit{.}

\bigskip

\noindent\textbf{Proof.} Starting from $\varepsilon^{(2)}:E^{(2)}%
\longrightarrow\mathrm{M}_{2}(K[z_{1},z_{2}]/(z_{1}^{2},z_{2}^{2}))$ and using%
\[
(K[z_{1},\ldots,z_{m}]/(z_{1}^{2},\ldots,z_{m}^{2}))[z]/(z^{2})\cong
K[z_{1},\ldots,z_{m},z_{m+1}]/(z_{1}^{2},\ldots,z_{m}^{2},z_{m+1}^{2}),
\]
iteration of the construction in Theorem 3.1 gives the desired
CT-representation. $\square$

\bigskip

\noindent\textbf{3.3.Proposition.}\textit{ If }$\varepsilon:E^{(m)}%
\longrightarrow\mathrm{M}_{t}(\Omega)$\textit{ is a (not necessarily CT-)
representation of }$E^{(m)}$\textit{ for some integers }$m,t\geq2$
\textit{over a commutative }$K$\textit{-algebra }$\Omega$\textit{, then
}$m\leq2t-1$\textit{.}

\bigskip

\noindent\textbf{Proof.} Let $S_{m}(x_{1},\ldots,x_{m})=\sum_{\pi
\in\mathrm{Sym}(m)}\mathrm{sgn}(\pi)x_{\pi(1)}\cdots x_{\pi(m)}$ be the
standard polynomial in $K\left\langle x_{1},\ldots,x_{m}\right\rangle $. The
well-known Amitsur--Levitzki theorem (see [R]) asserts that $S_{2t}=0$ is a
polynomial identity on $\mathrm{M}_{t}(\Omega)$. The existence of the
embedding $\varepsilon$ ensures that $S_{2t}=0$ is also an identity on
$E^{(m)}$. On the other hand%
\[
S_{m}(v_{1},\ldots,v_{m})=m!\;v_{1}\cdots v_{m}\neq0
\]
shows that $S_{m}=0$ is not a polynomial identity on $E^{(m)}$. If $2t\leq m$,
then $S_{m}=0$ follows from $S_{2t}=0$ (in any algebra). Thus we have
$2t\nleqslant m$. $\square$

\bigskip

\noindent\textbf{3.4.Theorem}("Cayley--Hamilton")\textbf{.} \textit{If }%
$A\in\mathrm{M}_{n}(E^{(m)})$\textit{ is an }$n\times n$\textit{ matrix, then
}$A$\textit{ satisfies an identity of the form}%
\[
A^{2^{m-1}n}+c_{1}A^{(2^{m-1}n)-1}+\cdots+c_{(2^{m-1}n)-1}A+c_{2^{m-1}n}%
I_{n}=0,
\]
\textit{where }$c_{i}\in K$\textit{, }$1\leq i\leq2^{m-1}n$\textit{.}

\bigskip

\noindent\textbf{Proof.} A combination of Theorem 2.1 and Corollary 3.2 gives
the identity. $\square$

\bigskip

\noindent\textbf{3.5.Remark. }In particular, we know the following about
embeddability of $E^{(3)}$ into matrix algebras. By Corollary 3.2, $E^{(3)}$
admits a CT-embedding into a $4\times4$ matrix algebra. On the other hand, it
has no CT-embedding into any $2\times2$ matrix algebra. The simple reason is
that for $c_{1},c_{2}\in K$, the $2\times2$\ Cayley--Hamilton identity%
\[
(v_{1}+v_{2}v_{3})^{2}+c_{1}(v_{1}+v_{2}v_{3})+c_{2}=0
\]
never holds in $E^{(3)}$, although such an identity should hold in $E^{(3)}$
by Theorem 2.1 if it CT-embeds into a $2\times2$ matrix algebra. We could not
decide whether $E^{(3)}$ has a CT-embedding into a $3\times3$ matrix algebra.

For $m=3$ and $t=2$ the inequality $m\leq2t-1$ in Proposition 3.3 becomes
equality, thus Proposition 3.3 does not contradict the existence of a non
CT-representation $\varepsilon:E^{(3)}\longrightarrow\mathrm{M}_{2}(\Omega)$.
Indeed, the assignments%
\[
1\longmapsto\left[
\begin{array}
[c]{cc}%
1 & 0\\
0 & 1
\end{array}
\right]  ,v_{1}\longmapsto\left[
\begin{array}
[c]{cc}%
0 & x\\
x & 0
\end{array}
\right]  ,v_{2}\longmapsto\left[
\begin{array}
[c]{cc}%
y & 2y\\
-2y & -y
\end{array}
\right]  ,v_{3}\longmapsto\left[
\begin{array}
[c]{cc}%
-2z & -z\\
z & 2z
\end{array}
\right]
\]
define a $K$-embedding $E^{(3)}\longrightarrow\mathrm{M}_{2}(K[x,y,z]/(x^{2}%
,y^{2},z^{2}))$. Since%
\[
v_{1}v_{2}v_{3}\longmapsto\left[
\begin{array}
[c]{cc}%
3xyz & 0\\
0 & 3xyz
\end{array}
\right]  ,
\]
this embedding is not a CT-representation.

\bigskip

\noindent\textbf{3.6.Remark.} Any matrix $A\in\mathrm{M}_{n}(E)$ (here
$\mathrm{char}(K)=0$ is essential) satisfies left and right Cayley--Hamilton
identities of the form%
\[
A^{n^{2}}+g_{1}A^{n^{2}-1}+\cdots+g_{n^{2}-1}A+g_{n^{2}}I_{n}=0,
\]%
\[
A^{n^{2}}+A^{n^{2}-1}h_{1}+\cdots+Ah_{n^{2}-1}+I_{n}h_{n^{2}}=0,
\]
and a central Cayley--Hamilton identity of the form%
\[
A^{2n^{2}}+u_{1}A^{2n^{2}-1}+\cdots+u_{2n^{2}-1}A+u_{2n^{2}}I_{n}=0,
\]
where $g_{i},h_{i}\in E,$ $1\leq i\leq n^{2}$ and $u_{j}\in E_{0},$ $1\leq
j\leq2n^{2}$ (see [Sz1]).

\bigskip

\noindent\textbf{3.7.Theorem.} \textit{The standard identity }$S_{2^{m}n}%
=0$\textit{ of degree }$2^{m}n$\textit{\ is a polynomial identity on
}$\mathrm{M}_{n}(E^{(m)})$\textit{.}

\bigskip

\noindent\textbf{Proof.} Since the natural extension%
\[
\left(  \varepsilon^{(m)}\right)  _{n}:\mathrm{M}_{n}(E^{(m)})\longrightarrow
\mathrm{M}_{n}(\mathrm{M}_{2^{m-1}}(K[z_{1},\ldots,z_{m}]/(z_{1}^{2}%
,\ldots,z_{m}^{2})))
\]
of $\varepsilon^{(m)}$\ (in Corollary 3.2) is a $K$-embedding and%
\[
\mathrm{M}_{n}(\mathrm{M}_{2^{m-1}}(K[z_{1},\ldots,z_{m}]/(z_{1}^{2}%
,\ldots,z_{m}^{2})))\cong\mathrm{M}_{2^{m-1}n}(K[z_{1},\ldots,z_{m}%
]/(z_{1}^{2},\ldots,z_{m}^{2}))
\]
satisfies $S_{2^{m}n}=0$ by the Amitsur--Levitzki theorem, the proof is
complete. $\square$

\bigskip

If $n=1$, then Theorems 3.4 and 3.7 are far from being sharp.

\bigskip

\noindent\textbf{3.8.Remark.} Using Theorem 5.5 of Domokos [D] and the fact
that $S_{m+2}=0$ is an identity on $E^{(m)}$, we obtain that $S_{(m+1)n^{2}%
+1}=0$ is an identity on $\mathrm{M}_{n}(E^{(m)})$.

\bigskip

The left regular representation of $E^{(m)}$ is an embedding%
\[
\lambda^{(m)}:E^{(m)}\longrightarrow\mathrm{End}_{K}(E^{(m)})\cong
\mathrm{M}_{2^{m}}(K)
\]
of $K$-algebras, where $\mathrm{End}_{K}(E^{(m)})$ is the algebra of all
$K$-linear maps of the $2^{m}$ dimensional vector space $_{K}E^{(m)}$ and
$\lambda^{(m)}(g):E^{(m)}\longrightarrow E^{(m)}$ is the left multiplication
by $g\in E^{(m)}$. The size of the matrix algebra is $2^{m-1}$ in Corollary
3.2, half of the size $2^{m}$ provided by the above left regular
representation. On the other hand, the base field $K$ is replaced by a much
bigger base ring ($K$-algebra) in Corollary 3.2. Since we cannot derive
Theorems 3.4 and 3.7 from the regular representation, our half sized embedding
is really better in many ways than the regular representation.

If we keep $K$ as the base field, then the next theorem gives a lower bound
for the matrix size of any possible embedding of $E^{(m)}$.

\bigskip

\noindent\textbf{3.9.Theorem.} \textit{If }$\lambda:E^{(m)}\longrightarrow
\mathrm{M}_{t}(K)$\textit{ is an embedding of }$K$\textit{-algebras, then
}$3\cdot2^{m-2}=2^{m-1}+2^{m-2}\leq\left\lfloor \frac{t^{2}}{4}\right\rfloor
+1$\textit{.}

\bigskip

\noindent\textbf{Proof.} Consider the $K$-subalgebra $E_{0}^{(m)}[v_{m}]$ of
$E^{(m)}$\ generated by the centre $E_{0}^{(m)}$ and the generator $v_{m}$.
Clearly, $E_{0}^{(m)}[v_{m}]$ is commutative, and any element of $E_{0}%
^{(m)}[v_{m}]$ is of the form $g_{0}+h_{0}v_{m}$, where $g_{0}\in E_{0}^{(m)}$
and $h_{0}\in E_{0}^{(m-1)}$. It follows that $\dim_{K}E_{0}^{(m)}%
[v_{m}]=2^{m-1}+2^{m-2}$. Since $\lambda(E_{0}^{(m)}[v_{m}])\subseteq
\mathrm{M}_{t}(K)$ is a commutative subalgebra of $\mathrm{M}_{t}(K)$, we can
apply Schur's inequality (see [M]):%
\[
2^{m-1}+2^{m-2}=\dim_{K}E_{0}^{(m)}[v_{m}]=\dim_{K}\lambda(E_{0}^{(m)}%
[v_{m}])\leq\left\lfloor \frac{t^{2}}{4}\right\rfloor +1.\text{ }\square
\]
\noindent\textbf{3.10.Remark.} Since M. Domokos and M. Zubor recently proved
(see [DZ]) that $\dim_{K}C\leq3\cdot2^{2r-2}$ for any commutative
$K$-subalgebra $C\subseteq E^{(2r)}$ of $E^{(2r)}$, the argument in the above
proof cannot be improved if $m=2r$ is even.

\bigskip

\noindent4. SKEW\ POLYNOMIAL\ RINGS AND EMBEDDINGS

\bigskip

For a ring ($K$-algebra) endomorphism $\sigma:R\longrightarrow R$, consider
the skew polynomial ring ($K$-algebra) $R[w,\sigma]$ in the skew indeterminate
$w$. The elements of $R[w,\sigma]$ are left polynomials of the form
$f(w)=a_{0}+a_{1}w+\cdots+a_{k}w^{k}$ with $a_{0},a_{1},\ldots,a_{k}\in R$.
Besides the obvious addition, we have the following multiplication rule in
$R[w,\sigma]$: $wr=\sigma(r)w$ for all $r\in R$ and
\[
(a_{0}+a_{1}w+\cdots+a_{k}w^{k})(b_{0}+b_{1}w+\cdots+b_{l}w^{l})=c_{0}%
+c_{1}w+\cdots+c_{k+l}w^{k+l},
\]
where
\[
c_{m}=\underset{i+j=m,i\geq0,j\geq0}{\sum}a_{i}\sigma^{i}(b_{j}).
\]
If $\sigma$ is an involution, then $w^{2}$ is a central element of
$R[w,\sigma]$: we have $\sigma(\sigma(r))=r$ and $w^{2}r=w\sigma
(r)w=\sigma(\sigma(r))w^{2}=rw^{2}$ for all $r\in R$, moreover $w^{2}$
commutes with the powers of $w$. Thus the ideal $(w^{2})$ of $R[w,\sigma]$
generated by $w^{2}$ can be written as $(w^{2})=R[w,\sigma]w^{2}%
=w^{2}R[w,\sigma]$. Consider the factor ring ($K$-algebra) $R[w,\sigma
]/(w^{2})$, then for any element $f(w)\in R[w,\sigma]$ there exists exactly
one left polynomial of the form $r+sw\in R[w,\sigma]$ in the residue class
$f(w)+(w^{2})$. Hence the elements of $R[w,\sigma]/(w^{2})$ can be represented
by linear left polynomials with coefficients in $R$, and the multiplication in
$R[w,\sigma]/(w^{2})$ is the following:
\[
(r+sw)(p+qw)=rp+(rq+s\sigma(p))w,
\]
where $r,s,p,q\in R$.

If $R=E=E_{0}\oplus E_{1}$ (or $R=E^{(m)}=E_{0}^{(m)}\oplus E_{1}^{(m)}$) is
the natural $\mathbb{Z}_{2}$-grading, then $\sigma(g)=\sigma(g_{0}%
+g_{1})=g_{0}-g_{1}$ defines a natural involution (here $g=g_{0}+g_{1}$ is the
unique presentation as a sum of an even and an odd element). It is easy to see
that $E[w,\sigma]/(w^{2})\cong E$ and $E^{(m)}[w,\sigma]/(w^{2})\cong
E^{(m+1)}$ as $K$-algebras.

We note that the idea of considering $R[w,\sigma]/(w^{2})$ comes from [SSz].

\bigskip

\noindent\textbf{4.1.Theorem}("Fundamental Embedding")\textbf{.}\textit{ For
an involution }$\sigma:R\longrightarrow R$\textit{, putting}%
\[
\mu(r+sw+(w^{2}))=\left[
\begin{array}
[c]{cc}%
r+(z^{2}) & sz+(z^{2})\\
\sigma(s)z+(z^{2}) & \sigma(r)+(z^{2})
\end{array}
\right]
\]
\textit{(with} $r,s\in R$\textit{) gives a }$K$\textit{-embedding }%
$\mu:R[w,\sigma]/(w^{2})\longrightarrow\mathrm{M}_{2}(R[z]/(z^{2}))$\textit{.}

\bigskip

\noindent\textbf{Proof.} We only have to prove the multiplicative property of
$\mu$:%
\[
\mu\left(  (r+sw+(w^{2}))(p+qw+(w^{2}))\right)  =\mu\left(  rp+(rq+s\sigma
(p))w+(w^{2})\right)  =
\]%
\[
\left[
\begin{array}
[c]{cc}%
rp+(z^{2}) & (rq+s\sigma(p))z+(z^{2})\\
\sigma(rq+s\sigma(p))z+(z^{2}) & \sigma(rp)+(z^{2})
\end{array}
\right]  =
\]%
\[
\left[
\begin{array}
[c]{cc}%
r+(z^{2}) & sz+(z^{2})\\
\sigma(s)z+(z^{2}) & \sigma(r)+(z^{2})
\end{array}
\right]  \cdot\left[
\begin{array}
[c]{cc}%
p+(z^{2}) & qz+(z^{2})\\
\sigma(q)z+(z^{2}) & \sigma(p)+(z^{2})
\end{array}
\right]  =
\]%
\[
\mu(r+sw+(w^{2}))\mu(p+qw+(w^{2})).\text{ }\square
\]

The following is a broad generalization of Theorem 3.1.

\bigskip

\noindent\textbf{4.2.Theorem.} \textit{Let }$\varepsilon:R\longrightarrow
\mathrm{M}_{t}(\Omega)$ \textit{be a CT-representation of }$R$\textit{ over a
commutative }$K$\textit{-algebra }$\Omega$\textit{ for some integer }$t\geq
2$\textit{. If }$\sigma:R\longrightarrow R$\textit{ is an involution, then
there exists an induced CT-representation }$\varepsilon^{\ast}:R[w,\sigma
]/(w^{2})\longrightarrow\mathrm{M}_{2t}(\Omega\lbrack z]/(z^{2}))$\textit{ of
the factor }$R[w,\sigma]/(w^{2})$\textit{.}

\bigskip

\noindent\textbf{Proof.} Consider the natural extension%
\[
\overline{\varepsilon}_{2}:\mathrm{M}_{2}(R[z]/(z^{2}))\longrightarrow
\mathrm{M}_{2}((\mathrm{M}_{t}(\Omega)[z])/(z^{2}))\cong\mathrm{M}_{2t}%
(\Omega\lbrack z]/(z^{2}))
\]
of $\overline{\varepsilon}:R[z]/(z^{2})\longrightarrow(\mathrm{M}_{t}%
(\Omega)[z])/(z^{2})$, where $\overline{\varepsilon}(r+sz+(z^{2}%
))=\varepsilon(r)+\varepsilon(s)z+(z^{2})$ for $r,s\in R$. Using the map
$\mu:R[w,\sigma]/(w^{2})\longrightarrow\mathrm{M}_{2}(R[z]/(z^{2}))$ in
Theorem 4.1, the composition $\varepsilon^{\ast}=\overline{\varepsilon}%
_{2}\circ\mu$ gives the induced CT-representation. Indeed, $\mathrm{tr}%
(\varepsilon^{\ast}(r+sw+(w^{2})))=(\mathrm{tr}(\varepsilon(r))+(z^{2}%
))+(\mathrm{tr}(\varepsilon(\sigma(r)))+(z^{2}))\in K+(z^{2}).$ $\square$

\bigskip

Let $\sigma:R\longrightarrow R$ be a ($K$-algebra) endomorphism with
$\sigma^{t}=1$ (such a $\sigma$\ is an automorphism). Now $w^{t}$ is a central
element of $R[w,\sigma]$ and the ideal $(w^{t})$ of $R[w,\sigma]$ can be
written as $(w^{t})=R[w,\sigma]w^{t}=w^{t}R[w,\sigma]$. We close this section
by mentioning the following generalization of Theorem 4.1 (see [Sz2]), which
seems to have applications in Galois theory.

\bigskip

\noindent\textbf{4.3.Theorem.}\textit{ For an endomorphism }$\sigma
:R\longrightarrow R$\textit{ with }$\sigma^{t}=1$\textit{,} \textit{putting}%
\[
\mu(r_{0}+r_{1}w+\cdots+r_{t-1}w^{t-1}+(w^{t}))=\left[  \sigma^{i-1}%
(r_{j-i})z^{j-i}+(z^{t})\right]  _{t\times t}%
\]
\textit{gives an embedding }$\mu:R[w,\sigma]/(w^{t})\longrightarrow
\mathrm{M}_{t}(R[z]/(z^{t}))$\textit{,} \textit{where the difference }
$j-~i\in\{0,1,\ldots,t-1\}$\textit{ is taken modulo }$t$\textit{, and the
element }$\sigma^{i-1}(r_{j-i})z^{j-i}+(z^{t})$\textit{ of the factor}
$R[z]/(z^{t})$\textit{ is in the }$(i,j)$\textit{ position of the }$t\times
t$\textit{ matrix }$\left[  \sigma^{i-1}(r_{j-i})z^{j-i}+(z^{t})\right]
_{t\times t}$\textit{. The trace of }$\left[  \sigma^{i-1}(r_{j-i}%
)z^{j-i}+(z^{t})\right]  _{t\times t}$ \textit{is in }$R^{\sigma}+(z^{t}%
)$\textit{, where }$R^{\sigma}=\{r\in R\mid\sigma(r)=r\}$\textit{ is the fixed
ring of }$\sigma$\textit{.}

\bigskip

\noindent REFERENCES

\noindent\lbrack D] M. Domokos, Eulerian polynomial identities and algebras
satisfying a standard identity, \textit{J. Algebra} 169(3) (1994), 913-928.

\noindent\lbrack DZ] M. Domokos and M. Zubor, Commutative subalgebras of the
Grassmann algebra, arXiv:1403.2916

\noindent\lbrack K] A. R. Kemer,\textit{\ }Ideals of Identities of Associative
Algebras, \textit{Translations of Math. Monographs}, Vol. 87 (1991), AMS,
Providence, Rhode Island.

\noindent\lbrack M] M. Mirzakhani, A simple proof of a theorem of Schur,
\textit{Amer. Math. Monthly} 105(3) (1998), 260-262.

\noindent\lbrack R] L.H. Rowen, Polynomial identities in ring theory\textit{,}
\textit{Academic Press}, New York, 1980.

\noindent\lbrack SSz] S. Sehgal and J. Szigeti, Matrices over centrally
$\mathbb{Z}_{2}$-graded rings, \textit{Beitr. Algebra Geom.} 43(2) (2002), 399-406.

\noindent\lbrack Sz1] J. Szigeti, New determinants and the Cayley--Hamilton
theorem for matrices over Lie nilpotent rings, \textit{Proc. Amer. Math. Soc.}
125(8) (1997), 2245-2254.

\noindent\lbrack Sz2] J. Szigeti,\ Embedding truncated skew polynomial rings
into matrix rings and embedding of a ring into 2x2 supermatrices, arXiv:1307.1783

\noindent\lbrack SzvW] J. Szigeti and L. van Wyk, Determinants for $n\times n$
matrices and the symmetric Newton formula in the $3\times3$ case,
\textit{Linear and Multilinear Algebra}, Vol. 62, No.8 (2014), 1076-1090.

\end{document}